\documentclass[12pt]{article}

\usepackage{epigraph,hyperref,fullpage}

\begin{document}

\title{Emotional Labor in Mathematics:\\ Reflections on Mathematical Communities,\\ Mentoring Structures, and EDGE}
\author{Gizem Karaali}
\date{\small Department of Mathematics, Pomona College, Claremont, CA, USA\\{\tt gizem.karaali@pomona.edu}}
%
%
\maketitle

\abstract{Terms such as ``affective labor'' and ``emotional labor'' pepper feminist critiques of the workplace. Though there are theoretical nuances between the two phrases, both kinds of labor involve the management of emotions; some acts associated with these constructs involve caring, listening, comforting, reassuring, and smiling. In this article I explore the different ways academic mathematicians are called to provide emotional labor in the discipline, thereby illuminating a rarely visible component of a mathematical life in the academy. 
Underlying this work is my contention that a conceptualization of labor involved in managing emotions is of value to the project of understanding the character, values, and boundaries of such a life. In order to investigate the various dimensions of emotional labor in the context of academic mathematics, I extend the basic framework of Morris and Feldman \cite{Morris} and then apply this extended framework to the mathematical sciences. Other researchers have mainly focused on the negative effects of emotional labor on a laborer's physical, emotional, and mental health, and several examples in this article align with this framing. However, at the end of the article, I argue that mathematical communities and mentoring structures such as EDGE help diminish some of the negative aspects of emotional labor while also accentuating the positives.}


\epigraph{The Bureau of Labor Statistics puts both ``diplomat" and ``mathematician" in the ``professional" category, yet the emotional labor of a diplomat is crucial to his work whereas that of a mathematician is not.
}{Arlie Russell Hochschild\\{\it The Managed Heart} \cite{TMH}, page 148.}

\section{Introduction}
\label{sec:1}

\parskip=7pt
\parindent=0pt

Sociologist Arlie Russell Hochschild launched the term ``emotional labor" into the mainstream with her 1983 book {\it The Managed Heart:\ Commercialization of Human Feeling}.  Inspired by Hochschild, and the large literature on emotional labor following her seminal work (see, for instance, \cite{Steinberg} for a careful review of this work up to the end of the twentieth century, or \cite{Grande} for a more recent volume), in this paper I will define emotional labor as any labor that involves the management of emotions, of the self or of others. For some nuances that might help the reader engage with contemporary literature on this theme, see Section \ref{S:ELandAL}. 

No matter how one defines emotional labor, according to Hochschild, the job of a mathematician is quite independent of this kind of labor (as per the epigraphed quote). Indeed many mathematicians enter this profession with similar illusions. For some it is even an appealing factor that the human contact required in many other professions would not be relevant. For others, emotional labor is a non-issue; they at some point decide they love and cannot live without mathematics, or they decide math is something they can do well enough to feed themselves and their loved ones; in either case, the emotional dimensions of labor do not come into play in their internal negotiations about future career plans. 

However many academic mathematicians soon find that their job entails emotional labor even if it is not part of the explicit job description.\footnote{ Here for reference is the official job description the United States government provides: mathematicians ``conduct research in fundamental mathematics or in application of mathematical techniques to science, management, and other fields", see \url{https://www.bls.gov/oes/current/oes152021.htm}, last accessed on February 25, 2019.} In this article I investigate the different ways mathematicians are called upon to provide emotional labor in the discipline. 
Even though many scholars have explored the features of emotional labor in academia, literature does not engage with the specific context and experiences of those in the mathematical sciences. In this article I probe the construct of emotional labor in the context of academic mathematics in order to shed light on this oft-neglected dimension of work in the discipline, and to highlight some aspects of it that might otherwise be missed. Thus this article may be viewed as a case study of sorts for emotional labor scholarship on the one hand and a reflection exercise for mathematicians on the other. 

In Section \ref{S:ELandAL}, I review the literature on emotional and affective labor, and building on prior work, I tease out a nine-dimensional framework that undergirds the discussion in the rest of the paper. In Section \ref{S:ELinA}, I apply this framework to the three contexts of teaching, academic service, and academic research. In this section I also begin to answer the related question: ``Who hears the call to emotional labor?" In Section \ref{S:EALinM}, I focus on the mathematical context and identify the kinds of emotional labor mathematicians are called to do, adapting the framework of Section \ref{S:ELandAL} and the examples of Section \ref{S:ELinA} into the mathematical sciences, and supplementing them with new ideas as needed.  In Section \ref{S:ELAinM-benefits}, I zero in on a specific category of work demanding emotional labor. Thus in this section I explore whether, and if so how, mathematical communities and mentoring structures (such as EDGE) may help diminish some of the negative aspects of emotional labor while also accentuating the positives. I also address more fully the ``who" question posed earlier. 
Section \ref{S:Final} wraps up this article, with a few brief remarks.



\vspace{-10pt}

\section{A Framework for Emotional Labor}
\label{S:ELandAL}

Theorists of labor have explored various types of work through the last century. Manual / physical labor, that is, work that engages the human body in the production of  commodities, has been central to most labor movements of the twentieth century. Intellectual / cognitive labor is of interest to many living in today's ``post-industrial" ``knowledge economy" environment. This paper is about a third type of labor, which as Hochschild put it bluntly, has been ``seldom recognized by those who tell us what labor is" \cite[page 197]{TMH}, also see \cite{Harding}. Indeed emotional labor is often viewed to be feminine and thus ``less-than".

In this paper we take any labor that involves managing the emotions of the worker or of those they interact with to be emotional labor. Some acts associated with this kind of labor include, but are not limited to, caring, listening, comforting, reassuring, and smiling. However there are two distinct components here: the self-management component remains internal, while the outward management of the emotions of the other (the client, the patient, the passenger, or, in the classroom setting, the student) is often more explicitly delineated and externally monitored by the employer. Even though both kinds of emotion-related work were labeled emotional labor by Hochschild in her \cite{TMH}, today these two are typically analyzed under different terms. 

In most contemporary scholarship, work that entails the monitoring and managing of the emotions of the laborer is called ``emotional labor"; see for example, \cite{Morris}, as well as the many articles in \cite{SteinbergFigartVolume}. Work that entails the creation or management of emotions in a designated other (or a designated group of others) is called ``affective labor", after \cite{Hardt}. These terms and the scholarship engaging with them are not uncontroversial; see for instance \cite{Schultz} for a critique of how the phrase  ``affective labor" might be used to create a gendered hierarchy of labor. Nonetheless, the study of emotions in the workplace in general has led to productive ways of thinking about the wellbeing of workers and more broadly about organizational behavior \cite{Fisher}. Furthermore, understanding how mental, manual, and emotional labor come together in today's work may lead to ``a potentially more comprehensive understanding of nature and social life" \cite{Harding}. Thus, it is my contention that a conceptualization of labor involved in managing emotions is of value to the project of understanding the character, values, and boundaries of a mathematical life in the academy. In the following, therefore, I will use the term ``emotional labor" to capture both types of emotion work, pointing out explicitly the distinct aspects of different types when needed. 


A formal framework to conceptualize the type of emotional labor involving the management of the emotions of the laborer is presented in \cite{Morris}, where four dimensions are proposed:
\begin{enumerate}
\item[(a)\textsubscript{self}] frequency of appropriate emotional display, 
\item[(b)\textsubscript{self}] attentiveness to required display rules, 
\item[(c)\textsubscript{self}] variety of emotions required to be displayed, and 
\item[(d)\textsubscript{self}] emotional dissonance generated as the result of having to express organizationally desired emotions not genuinely felt.
\end{enumerate}
The authors then argue that ``although some dimensions of emotional labor (e.g., variety of emotions that are
displayed) are likely to be associated with higher emotional exhaustion, it is mainly emotional dissonance that is likely to lead to lower job satisfaction." Some of what follows will have resonances with this perspective; see in particular Section \ref{S:ELAinM-benefits}.

A comprehensive framework for an analysis of emotional labor in the context of academic mathematics will also need to account for the dimensions of the type of emotion work that involves the creation and management of desired emotions in designated others. Analogous to the four dimensions above, I propose the following:
\begin{enumerate}
\item[(a)\textsubscript{others}] frequency of instances of management of the emotions of designated others, 
\item[(b)\textsubscript{others}] attentiveness to designated others' current emotions, 
\item[(c)\textsubscript{others}] variety of emotions one is required to engender or sustain in designated others, and 
\item[(d)\textsubscript{others}] emotional burden generated as the result of having to focus on designated others' emotions at the expense of other priorities or personal values.
\end{enumerate}

To the above eight dimensions, I will add a ninth that does not require the management of displayed emotion and yet is very much related to internal self-directed emotion management:
\begin{enumerate}
\item[(e)\textsubscript{self}] internal self-management of emotion required to continue to perform effectively in the job.
\end{enumerate}
This dimension of emotion work, related in various ways to the mental health of the laborer, is typically not monitored by the employer and yet is absolutely crucial to the employee's performance and sustained effectiveness. See \cite{Wharton} for a review of the psychosocial consequences of the self-monitoring of emotions, and \cite{Pugli} for recent work that points toward the various significant effects of self-focused emotion management on the mental health and well-being of workers. 

\vspace{-10pt}

\section{Emotional Labor in the Academy}
\label{S:ELinA}

There are many ways of doing emotional labor in the academy. See \cite{Taylor} for a ``contemporary account of what it means to experience and feel academia, as a privilege, risk, entitlement, or failure" [page 1]. Bellas \cite{Bellas} offers a critique of the rewards system within the academy that values the ``masculine" aspects of the job (research and administration) over the ``feminine" aspects (teaching and service). However she also points out that emotional labor plays a significant role in all these areas of academic work. 

In the three subsections of this section I explore some examples of the emotional labor involved in the trifecta of teaching, research, and service in terms of the nine dimensions described above. This will be a general exploration; I will leave the analysis of the specific context of the academic mathematician to the next section.

In the examples of this section (and the next), I intentionally focus almost exclusively on the burdens that may accrue from the emotional dimensions of teaching, service, and research. I do not deny that academics also get a lot of emotional satisfaction from their work, and that for many, some of these aspects serve as the highlight of their careers and sustain them for years. In fact in Section \ref{S:ELAinM-benefits}, I will come back to teaching and service in particular, and zero in on some of the positives of the emotional labor involved in these two aspects of an academic career. 

\vspace{-10pt}

\subsection{Emotional dimensions of teaching}
\label{SS:ETeaching}

Teaching involves a significant amount of emotional work. Sarah Rose Cavanagh in \cite{Cav} explores in depth the emotional dimensions of teaching and learning; see especially pages 102--108 where the emotional labor of teaching is explicitly discussed. N\"{a}ring and coauthors in \cite{Naring} explore specifically the relationship between the emotion work of teachers and their emotional exhaustion. For this paper, however, it will be sufficient to specifically point out examples of emotional labor related to teaching that can be described in terms of the nine dimensions of Section \ref{S:ELandAL}. 

During the academic year, professors typically meet their students in the classroom a few times a week ((a)\textsubscript{self}). During these regular sessions, they often aim to display effortless expertise, enthusiasm for the discipline, and joy of teaching ((c)\textsubscript{self}). These impressions are not always easy to sustain, and for those instructors who are traditionally underrepresented in the professoriate, they may be somewhat difficult to sustain simultaneously. (``She is so bubbly enthusiastic about her topic! She must not really know what she is talking about.") Thus the professor must often pay close attention to carefully balancing the displayed emotions ((b)\textsubscript{self}). This balancing act is often difficult. Too much enthusiasm might be counterproductive. One must be perceived as professional and yet friendly, charismatic and yet not too distant, and so on. This delicate performance aspect of teaching might be additionally difficult for the typical introvert academic, who might be exhausted by putting on a show for the students every other day. 

During office hours, professors regularly interact with students and need to manage their feelings ((a)\textsubscript{others}). They need to attend to feelings of helplessness, distrust, and apathy, and find ways of supporting students' confidence, interest, and enthusiasm ((b)\textsubscript{others},(c)\textsubscript{others}). They may need to turn on this others'-feelings focus, day after day, even when they themselves are not feeling emotionally healthy or when they have other needs of their own ((d)\textsubscript{others}). Beyond the standard requirements of office hour performance, occasionally professors find themselves in the roles of therapist, mother, experienced older brother, or wise elder, where they are expected to help students process and manage their emotions about various life issues. 

\parskip=9pt

In academia a disproportionate amount of teaching-related emotional work falls on the shoulders of female professors; in particular \cite{Alay} highlights the extra burden on female faculty. Add to this the gender bias in student evaluations of teaching, which have been explored at least since the publication of \cite{Basow}. Such extra work demands and special requests for favors as well as concerns about student evaluations certainly contribute to the burden of emotional work related to teaching, in particular in the dimensions of (b)\textsubscript{self}, (d)\textsubscript{self}, (e)\textsubscript{self}, (b)\textsubscript{others}, (d)\textsubscript{others}.

{It should be noted at this point that several statements made in this paper about the added burden of gender on female faculty in general and female mathematicians in particular may apply also to other minoritized or underrepresented groups. Most of the scholarship I have come across on the topic of emotional labor works with gender as the main variable. Scholarship that works with other variables such as race, ethnicity, and socioeconomic status in the context of the academy does not always engage explicitly with terms such as ``emotional labor" or ``affective labor". This of course does not mean that there is no extra burden on other minoritized groups in the academy. See for instance \cite{Thomp} for a discussion of the types of added stress African American faculty live with, and \cite{Lechuga} for an investigation into how underrepresented faculty may be managing their emotions. What is called ``emotional drain" in \cite{SVT} experienced by women of color in the academy is also very much in the same territory.}

\vspace{-10pt}

\subsection{Emotional dimensions of academic service}
\label{SS:EService}

There is evidence that supports the claim that women do significantly more academic service than men; research also seems to imply that the main source of discrepancy is in the amount of internal service \cite{Guarino}. Furthermore Bellas \cite{Bellas} argues that there may be ``greater demands for emotional labor on women as women attempt to convey, justify, and legitimize their contributions and, indeed, their presence." Noting how women's ways of leading and communication are often devalued, she concludes her section on emotional labor and academic service with ``All this may place additional demands for emotional labor on faculty women, especially women in the lower ranks and women of color, who may suffer devaluation in interactive as well as financial contexts." Once again readers may assume that the above will apply to other minoritized groups.

Two significant components of internal service are student advising and committee work \cite{Bellas}. Emotional labor related to student advising is much akin to the emotion work a professor does during office hours, even though often the others'-feelings focus is even more dominant in advising. Advising professors must remain in tune with their students' emotional needs and general emotional condition in any advising session ((b)\textsubscript{others}). Furthermore, they should preserve their ``caring professional" presentation ((b)\textsubscript{self}). 

In faculty committees, which involve work on the emotions of both self and others, professors often meet regularly to discuss matters that are either too small or too large to be resolved by said committees. The frequency of these meetings ((a)\textsubscript{self},(a)\textsubscript{others}), especially if in inverse proportion to their actual effectiveness, may lead to burnout and sometimes apathy. However the main emotional dimensions of this kind of academic work involve interpersonal relations between colleagues.  One should, for instance, make sure to help others feel good about themselves or at least not offend their sensibilities too much ((b)\textsubscript{others},(c)\textsubscript{others}). One should also attend to seeming interested and competent ((b)\textsubscript{self},(c)\textsubscript{self}). Faculty can feel emotionally burnt out if they find themselves assigned to committees whose work does not interest or challenge them, or if they find themselves on committees where their contributions are not valued. If they want to be good team players, they still feel the pressure to seem interested or at least act as if they care, all the while not feeling that way at all ((d)\textsubscript{self}). If on top of all this, their fellow committee members are high-maintenance folks who need extra emotion management, those professors who feel obliged to perform said management may be additionally burdened ((d)\textsubscript{others}). 

Undeniably some faculty find academic service outlets that are emotionally very fulfilling for them. Here, as mentioned in the preamble to Section \ref{S:ELinA}, I focused exclusively on the negative dimensions of emotion work related to service. See Section \ref{S:ELAinM-benefits} for an exploration
of the positive features to complement this discussion. 

%
%
%

\vspace{-10pt}

\subsection{Emotional dimensions of academic research}
\label{SS:EResearch}

Bellas \cite{Bellas} focuses on how emotional work relates to academic research in a section titled ``Emotional Labor and Research". In the following I mainly follow her examples, coding them according to the framework described earlier in Section \ref{S:ELandAL}. 

Certain types of research, in particular social science research which involves issues of personal relevance to the researcher, challenge the researcher to remain neutral and objective, or at least conscious of their biases; this might be emotionally challenging ((e)\textsubscript{self}). This, together with the expectation of neutrality in the presentation of the final product of the work, might lead to emotional dissonance ((d)\textsubscript{self}). In quantitative or empirical research work, one might still find emotional labor lurking in the background. 
If research involves interview or experiment participants, then the researcher may need to manage the emotions of  said participants ((b)\textsubscript{others}-(d)\textsubscript{others}). In all work, researchers need to remain vigilant against wishful thinking and overly optimistic interpretations of experimental results and other data ((e)\textsubscript{self}). 
If the research involves other researchers, such as training assistants, then the researcher once again will need to manage emotions of others, and this time probably at regular intervals ((a)\textsubscript{others}-(d)\textsubscript{others}).

Researchers need to attend to how the audience may view them during a conference talk, making sure to present a professional, competent, and yet interesting persona, though they may not really feel that way ((b)\textsubscript{self}-(c)\textsubscript{self}).
When submitting work for publication or proposals for grants or presentations, they should present themselves as competent and confident ((b)\textsubscript{self}-(c)\textsubscript{self}). If said work is rejected in a dismissive or rude manner, they must  pretend to be mature and generous and graciously take the given feedback 
((d)\textsubscript{self}). This type of effort might also be needed to handle certain audience members during conference presentations. 

\parskip=7pt

\vspace{-10pt}

\section{Emotional Labor in Mathematics}
\label{S:EALinM}

When I became a mathematician, I knew that my job would not involve much manual / physical labor.\footnote{ This is not to claim that mathematicians are disembodied workers. When I had a minor shoulder injury and had visions of not being able to use the chalk board for several weeks, or during that stressful time when I lost my voice unexpectedly, I very clearly noted the physical aspects of my role in the academy.} If I did think in terms of labor economics at all, I simply assumed that I would be part of today's knowledge economy, where I would be contributing to the production and dissemination of mathematical knowledge. But what was not obvious to my naive young self is that both research and education are a part of what is called the ``service-providing sector''; see for instance the classification offered by the United States Bureau of Labor Statistics \cite{BLS}. And this sector is today, possibly even more so than it had been during the writing of \cite{TMH}, the largest sector that demands emotional labor from its participants. 

In Sections \ref{SS:ETeaching}-\ref{SS:EResearch}, I explored specific ways in which the three main parts of an academic career (teaching, research, and service) might involve emotional labor. In this section I focus more explicitly on the context of the mathematical sciences. Note that in this section, too, my emphasis will be 
almost exclusively on the negative burden of these aspects of academic work. 

Perhaps it is natural that teaching mathematics is intrinsically emotional \cite{Boylan, Williams}. Students come into our classrooms with many emotions about mathematics. Some have math anxiety, some have self-doubt, some have a confidence level which may not serve them well in their next course. Students also bring along non-mathematical emotions, which contribute in all sorts of ways to how they engage with our content and pedagogy. If they just broke up with a partner, or if they have a sick relative, or if they are anxious about paying the next month's rent, their classroom participation as well as their learning will be impacted. Thus it makes sense that emotional labor, in particular others'-feelings oriented work along the dimensions of ((b)\textsubscript{others}-(d)\textsubscript{others}), makes up a significant portion of the work of a teaching mathematician.

Most mathematics professors can identify several familiar aspects of teaching described in Section \ref{SS:ETeaching}, if not in their own experiences, then in some of their colleagues'. In particular many departments have that one professor whose office hours tend to turn into what seem like therapy sessions from outside the door. This professor, more often than not, belongs to a minoritized group in the discipline; it can be a woman or a person of color for instance. It is clear that the amount of teaching-related emotional labor an individual professor takes on is not independent of the identity of that said professor. In particular, I have already mentioned in Section \ref{SS:ETeaching} that the distribution of emotional labor related to teaching does not seem to be gender-neutral;\linebreak there is no reason to expect that the situation will be different in mathematics.\footnote{ I am by no means suggesting that this is a desirable situation; nor am I asserting that this is a choice made by individual instructors. It is often the case that certain faculty find themselves in these situations. It is my belief that nobody should be forced to do more emotional labor than they are willing to do. Unfortunately, faculty from minoritized groups often face the dilemma of either doing the extra emotional work and not being respected for it or rejecting doing the extra emotional work and then suffering the consequences of that decision.}

The growing focus on pedagogies that emphasize student voice, student activity, and student agency in the mathematics classroom rather than a charismatic instructor's perfect presentation might be a reflection of the emotional labor of non-dominant groups in this area. Indeed professors from non-dominant groups may find that the ``traditional ways of being a mathematics professor" do not work for them. That is, the caricature of a genius mathematician staring at the chalk board or scanning the sea of nameless faces while delivering a flawless lecture may not be the ideal way for all professors to connect with and teach all students. Thus today's mathematics instructors, especially those from non-dominant groups, might tend toward teaching pedagogies that involve more others'-feelings focused emotional work, possibly thus lightening the self-directed emotional work load and the emotional dissonance that might accompany that kind of teaching. In \cite{Steurer}, Steurer describes how inquiry-based learning has resonated with her and allowed her to more naturally handle the emotion work of teaching. Conceiving of teaching mathematics as radical care \cite{SM} might be another way to reconcile these dissonances. 

Academic service of mathematicians is in many ways similar to that of other academics. Though some disciplines may be more open to mentoring and outreach activities than mathematics, and others might more naturally lead to service to the community in other ways (such as an engineering faculty member serving as a pro-bono consultant to a local water treatment facility), I believe that the features of emotion work described in Section \ref{SS:EService} capture many mathematicians' service experiences. Of course 
some may not have such uniformly negative experiences. There are many ways to serve, and some of these can be extremely fulfilling. For example Khadjavi and I worked hard for years to put together a collection of resources for mathematics instructors who want to incorporate social justice issues into their courses \cite{KK1, KK2}. This process was emotionally draining and yet also very satisfying. I also find journal editing work to be completely exhausting but also very fulfilling. I will come back to service in Section \ref{S:ELAinM-benefits} and there will address this point explicitly.

Thinking of academic research in mathematics, it is clear that mathematicians occasionally engage in research that is personally relevant to them. See for instance \cite{Kolba}, where Kolba describes research she undertook to better understand twin pregnancies, and \cite{Berger}, where Berger describes her research related to Down's Syndrome. Similarly the work I did with Glass on school districting (see for instance \cite{Glass}) also came out of personal investment in the topic. Nonetheless, I believe I need to move beyond the emotional dimensions already described in Section \ref{SS:EResearch} to fully engage with the emotional dimensions of mathematical research. 

To that end, I will refer to Weidman's list of the four emotional challenges of a mathematical life \cite{Weidman}:

\smallskip
\noindent{\it``First of all, the mathematician must be capable of total involvement in a specific problem."} That is, mathematics research often demands full focus for extended periods of time, and this is not only mentally exhausting but also emotionally draining. One might feel that one needs to withdraw from other interests, or else one is not doing enough. Each of the self-management dimensions ((a)\textsubscript{self}-(e)\textsubscript{self}) comes into play here.

\smallskip
\noindent{\it``Second, the mathematician must risk frustration. Most of the time, in fact, he finds himself, after weeks or months of ceaseless searching, with exactly nothing: no results, no ideas, no energy."} A lot of mathematics research work leads to no results of significance. Add to this the challenges of getting published once one does have a significant result, which, for the not-yet-thick-skinned, can get especially disorienting and discouraging. Similarly these are self-management focused ((b)\textsubscript{self}-(e)\textsubscript{self}). 

\smallskip
\noindent{\it``Next, even the most successful mathematician suffers from lack of appreciation."} The mathematics community proudly celebrates its geniuses, but celebrity and genius are fickle \cite{Karaali}. After all, what more can you do once you win a Fields medal? Anything after that will be a let-down. Even those who feel appreciated by their mathematical colleagues may suffer from a dearth of appreciation from family and friends, and the world outside the mathematical one might be totally immune to mathematical glory. This might lead to serious emotional strain ((e)\textsubscript{self}).

\smallskip
\noindent{\it``Finally, the mathematician must face the fact that he will almost certainly be dissatisfied with himself."} Somewhat a corollary of the above, this means that mathematics is huge and the contribution of each individual mathematician is just a small speck. Whatever one does will be small change when compared to some of the giants. Once again this type of dissatisfaction might lead to serious emotional strain ((e)\textsubscript{self}).

\smallskip

Weidman's four challenges all require mathematicians to manage their own emotions, corresponding to the dimensions ((a)\textsubscript{self}-(e)\textsubscript{self}). The mathematician who cannot handle all of them at least halfway successfully at least some of the time is bound to be miserable. One might reject some of these challenges as myths and try to disentangle oneself from their hold (see \cite{Harron} for a call for this kind of a rejection), but that too takes emotional work as these ideas and ideals are quite solidly built into the culture of the discipline.

So far I have explored some specific ways academic mathematicians might engage in emotional labor. Are there any other kinds of emotional labor mathematicians might be called to do? And just who gets to hear that call? What are the consequences of hearing that call? I focus on these questions in the next section. 

\vspace{-10pt}

\section{Mathematical Communities, Mentoring Structures, and EDGE}
\label{S:ELAinM-benefits}

The idea that emotional labor has various human costs is not new; neither is the idea that it is not uniformly a negative for the particular laborer engaged in it, see for instance \cite{Wharton}. In this section I explore some of the positive possibilities related to emotional labor in the context of an academic mathematical life. In particular I reflect upon mathematical communities, mentoring structures, and EDGE.

Since 2008, the American Mathematical Society has supported Mathematics Research Communities (MRC). The MRC is ``a professional development program offering early-career mathematicians a rich array of opportunities to develop collaboration skills, build a network focused in an active research domain, and receive mentoring from leaders in that area" \cite{AMSMRC}. In the last few years, the Association for Women in Mathematics has led or supported research networking conferences for women in various fields; see \cite{AWMadvance} for a list of research networks supported in various ways by the AWM. These programs, and others like them, are all spearheaded by mathematicians who feel called to do the work to create networks, connect people, mentor young mathematicians, and make our community a more welcoming and supportive place for more people. 

The Mathematical Association of America has two programs for mentoring junior mathematics faculty. The first one of these, Project NExT (New Experiences in Teaching), ``is a professional development program for new or recent Ph.D.s in the mathematical sciences" addressing ``all aspects of an academic career: improving the teaching and learning of mathematics, engaging in research and scholarship, finding exciting and interesting service opportunities, and participating in professional activities. It also provides the participants with a network of peers and mentors as they assume these responsibilities" \cite{NExT}.  The MAA Mentoring Network is another mentoring program ``aimed at connecting early career mathematicians with experienced mentors working in mathematics" \cite{MAAMentorNetwork}. 

There are many other mentoring structures built around academic mathematics. One might count among these:
\begin{enumerate}
\item The e-Mentoring Network, hosted by AMS blogs, available at \url{https://blogs.ams.org/mathmentoringnetwork/},
\item The Infinite Possibilities Conference, ``a national conference designed to promote, educate, encourage and support minority women interested in mathematics and statistics" \cite{IPC},
\item NSF Mathematics Institutes' Modern Math Workshop at SACNAS, ``a pre-con-ference workshop held at the SACNAS National Conference, intended to encourage undergraduates, graduate students and recent PhDs from underrepresented minority groups to pursue careers in the mathematical sciences and build research and mentoring networks"  \cite{MMW}. 
\end{enumerate}
And of course one cannot forget EDGE. ``The EDGE Program is administered by the Sylvia Bozeman and Rhonda Hughes EDGE Foundation with the goal of strengthening the ability of women students to successfully complete PhD programs in the mathematical sciences and place more women in visible leadership roles in the mathematics community. Along with the summer session, EDGE supports an annual conference, travel for research collaborations, travel to present research and other open-ended mentoring activities" \cite{EDGE}.

The work involved in each of these programs is varied, but there is a significant emotional labor component. The main job is to connect people to one another, and though a lot of the emotion work is distributed over a large number of people, the main program organizers do a large chunk of it.  There is much emotion work that involves the management of the emotions of others; in particular many of the junior mathematicians participating might be feeling insecure and lost or at least mildly confused. The emotional labor involved is mainly about making sure these participants feel a sense of belonging, and a sense of confidence and realistic optimism about their future in academic mathematics. 

Digging deeper, one can see the resonances with the types of emotional labor described in Sections \ref{SS:ETeaching}-\ref{SS:EService}. In particular the types of work involve teaching and service. However, people involved do the work willingly. They basically self-select into these roles. This is perhaps one of the main reasons why the emotional labor involved, though still highly burdensome in any objective sense of the word, does also help them feel nourished and fulfilled. 

However there are two other reasons I believe. 

First the people who put their time and energy into these programs feel called to do this type of work because they believe ideologically and philosophically that it is the right thing to do. Their political and ethical framing of the world puts them in the position to value this kind of work, and this in turn makes the work feel more endurable, more meaningful, and even more joyful.\footnote{ This certainly applies to other service work, such as my editing work with Khadjavi on \cite{KK1, KK2} and K-12 outreach activities of various mathematicians; see for instance \cite{WZ}. Similarly blogging and hosting other networking sites is a valuable service contribution, where most of the time people who do the needed work engage in it because of political and ethical goals. See in particular the e-Mentoring Network mentioned above as well as the inclusion/exclusion blog of the American Mathematical Society, available at \url{https://blogs.ams.org/inclusionexclusion/}.} ``Meaningful work" is a catchy phrase; see \cite{Chalofsky} for a working framework for it that revolves around three themes (sense of self, the work itself, and sense of balance), and see \cite{Steger} for more on the benefits of meaningful work for the laborer. But even leaving related scholarship aside, it is easy to understand how meaningful work can transform strenuous emotional labor into pleasurable and desirable labor. Thinking of the nine dimensions proposed in Section \ref{S:ELandAL}, one can see that the dimensions of emotional labor activated mainly involved are (a)\textsubscript{self}-(c)\textsubscript{self} and (a)\textsubscript{other}-(d)\textsubscript{other}. There is little self-deception or misrepresentation of feelings, and there is more or less no emotional dissonance. So perhaps the individuals are exhausted at the end of the day, but they sleep well.\linebreak  This resonates with the work of Morris and Feldman in \cite{Morris}, who found that emotional dissonance was the main component of emotional labor that led to job dissatisfaction; see the relevant quote in Section \ref{S:ELandAL}. 

Secondly and perhaps relatedly, there is often a shared identity component to the decision to dedicate time and energy to a program of this kind. This makes the work meaningful and the emotional dissonance minimal, yes, and in all these ways, this reason may seem similar to the first. But what makes this different is how it interacts with the others'-feelings focused labor dimensions (a)\textsubscript{other}-(d)\textsubscript{other}. The shared identity makes the emotional labor of managing the designated others' emotions a lot easier, as the individual has a better understanding of said emotions of those designated others. This of course does not mean that white women necessarily make the best mentors for white women, black men necessarily make the best mentors for black men, and gay Latinas necessarily make the best mentors for gay Latinas. But it is natural to expect that shared identity makes aspects of the involved emotional labor much easier. 

\vspace{-10pt}

\section{Concluding Thoughts}
\label{S:Final}

The EDGE Program is an example of projects that demand emotional labor but that also contribute significantly to the well-being of both its participants and its leaders and mentors. In fact several participants become mentors and leaders themselves; see \cite{Ale}. Though emotionally exhausting, these projects can continue because they fulfill several needs of those that work on them. Instead of expending their energy and emotional well-being on trying to run their departmental and campus committee meetings smoothly despite colleagues who expect them to make coffee or just smile and nod, and these mathematicians can put their efforts into projects that help them connect with each other, find meaning, and thrive in the academy. 

On a more pragmatic level, there is perhaps more that can be said. One of the anonymous reviewers of this paper wrote that ``the analysis [done here], which has intrinsic value in making visible labor that is generally hidden and unrewarded, has even more value-added as a framework to support these programs." I am sincerely energized by the possibility that this work may help these programs in one way or another. Indeed it may be of some value to the powers that be to comprehend that these programs benefit not only their participants but also their organizers, in ways that enrich their mathematical lives and naturally lighten some of the load of the emotional labor that comes with the territory. 

There are several directions where this work could be further enhanced. For example, as one of the reviewers mentioned, there is significant variation across institutions of higher education in the United States in terms of how the three pillars of academic work (teaching, service, research) are to be valued and prioritized. This might add significant complications to the way we can talk about  mathematicians' emotional labor. For instance, mathematics faculty in research-focused institutions might be concerned more about their research productivity and might have to perform more emotional work related to this aspect than others (in particular in the categories of ((b)\textsubscript{self}-(e)\textsubscript{self}) as described at the end of Section \ref{SS:EResearch}), while mathematicians in more teaching- and service-oriented institutions might find that their work involves more emotion work along the dimensions of ((b)\textsubscript{others}-(c)\textsubscript{others}) as described in Sections \ref{SS:ETeaching}-\ref{SS:EService}.

As another possible direction, consider that research suggests that ``foreign-born women faculty members' patterns of engagement in work activities contradict the gendered division of labor in academia" \cite{Mami}. Foreign-born faculty make up a significant percentage of the overall American academia (in 2009, nonresident aliens made up 11.5\% of the 11,599 new tenure-track faculty members at four-year institutions in the United States, while Asian Americans made up 10.5\%, African Americans made up 0.5\%, and Hispanics made 0.4\% \cite{Mami2}). Emotional labor in the context of foreign-born faculty might look different; think for instance about what kinds of different challenges such faculty might face in the contexts of teaching and advising, and how the benefits described in Section \ref{S:ELAinM-benefits} may apply to them in different ways. As many U.S.-based mathematicians are also foreign-born (cf.\ \cite{AIC, Lin}), parts of the discussion in this article on emotional labor in mathematics might need some modifying accordingly. 


\noindent
{\bf Acknowledgements.}
I would like to thank Sarah Bryant, Amy Buchmann, Susan D'Agostino, Michelle Craddock Guinn, and Leona Harris for giving me the opportunity to contribute to their edited volume. I also very much appreciated their open-mindedness about the uncategorizability of my project and their input and support through the writing process. I should also like to thank the two anonymous reviewers, who offered insightful feedback on the original version of this manuscript, and Brian Katz, who generously read two drafts of the paper and offered his constructive criticism. I believe that their input has much improved this paper. The imperfections that remain are of course all my responsibility; this is the best I could do at this time.

\end{document}